\documentclass{article}
\usepackage{graphicx}
\usepackage{fullpage}
\usepackage{latexsym}
\usepackage{amsmath}
\usepackage{amssymb}
\usepackage{bm}
\usepackage{wrapfig}
\usepackage[T1]{fontenc} 
\usepackage[english]{babel}
\usepackage{amsthm}
\usepackage{float}
\usepackage{subfigure}
\usepackage{fancybox}
\usepackage{fancyhdr}
\usepackage{epstopdf}
\numberwithin{equation}{section}

\newtheorem{theorem}{Theorem}
\numberwithin{theorem}{section}
\newtheorem{prof}{Proof}
\theoremstyle{definition}
\newtheorem{definition}{Definition}[section]
\newtheorem{corollary}{Corollary}[theorem]

\theoremstyle{remark}

\title{Associated curves in $E^3$ from a different point of view}
\author{S\"{u}leyman \uppercase{\c{S}enyurt}\\
Ordu University / Faculty of Arts and Science\\
Department of Mathematics\\
https://orcid.org/0000-0003-1097-5541
	\and
	Davut CANLI$^*$ \\
	Ordu University / Faculty of Arts and Science\\
Department of Mathematics\\
	https://orcid.org/ 0000-0003-0405-9969
		\and
	Kebire Hilal AYVACI \\
	Ordu University / Faculty of Arts and Science\\
Department of Mathematics\\
	https://orcid.org/0000-0002-5114-5475
	}

\date{}

\begin{document}

\maketitle%

\begin{abstract}
In this paper, tangent, principal normal and binormal wise associated curves are defined such that each of these vectors of any given curve lies on the osculating, normal and rectifying  plane of its mate, respectively. For each associated curves a new moving frame and the corresponding curvatures are found, and in addition to this the possible solutions for distance functions between the curve and its associated mate are discussed. In particular, it is seen that the involute curves belong to the family of tangent associated curves, the Bertrand and the Mannheim curves belong to the principal normal associated curves. Finally, as an application, we present some examples and map a given curve together with its mate and their frames.
\end{abstract}

\section{Introduction}

In differential geometry, curves are named as associated if there exist a mathematical relation among them.
Some of those known as involute-evolute curves, Bertrand curves, Mannheim curves and more recently the successor curves are the ones on which the researchers most studied \cite{bertrand}, \cite{menninger}, \cite{oneillb}, \cite{wangLiu}. For such curves, the association is based upon the Frenet elements of the curves. There have been other studies using different frames such as Darboux and Bishop to associate curves, as well (\cite{kayaonder1}, \cite{kayaonder2}, \cite{korpinar}, \cite{masal}, \cite{yilmazAykut}). 
From a distinct point of view, Choi and Kim (2012), introduced new associated curves of a given Frenet curve as the integral curves of vector fields. \cite{choiKim}.
\c{S}ahiner, on the other hand, established direction curves of "tangent" and "principal normal" indicatrix of any curve and provided some methods to portray helices and slant helices by using these curves in his studies, \cite{sahiner1}, \cite{sahiner2}, respectively.
In this study, we introduce another Frenet frame based associated curves such that the tangent, the principal normal and the binormal vectors of a given any curve lies on the osculating, normal and rectifying  plane of its partner, respectively. For each associated curves a new moving frame is established and the distances between the curve and its offset are given. In particular, it is seen that the involute curves belong to the family of tangent associated curves. In addition, some traces of the Bertrand and Mannheim curves are found while examining principal normal and binormal associated curves. Finally we provided some examples to illustrate the intuitive idea of this paper.

The Frenet vectors and curvatures of the curve, $\alpha$ together with the well known Frenet formulas are given as:
\begin{equation} \label{defns}
T(s)=\frac{\alpha '(s)}{\| \alpha '(s)\| },\,\,
N(s)=B(s)\wedge N(s),\,\,
B(s)=\frac{\alpha '(s)\wedge \alpha ''(s)}{\| \alpha '(s)\wedge \alpha ''(s)\| }, 
\end{equation}
\begin{equation}\label{curvatures}
\kappa(s) = \frac{\|\alpha'(s)\wedge \alpha''(s) \|}{\|\alpha'(s) \|^3}, \,\, \tau(s)= \frac{\langle \alpha'(s)\wedge \alpha''(s), \alpha'''(s)\rangle}{\|\alpha'(s)\wedge \alpha''(s) \|^2}.
\end{equation}
\begin{equation} \label{freneForms}
\vec{T^{{'}}}(s)=\kappa (s)\vec{N}(s),
\vec{N^{{'}}}(s)=-\kappa (s)\vec{T}(s)+\tau (s)\vec{B}(s),
\vec{B^{{'}}}(s)=-\tau (s)\vec{N}(s).
\end{equation}
where, $\vec{T}$, $\vec{N}$, $\vec{B}$, $\kappa$ and $\tau$ are called the tangent vector, the principal normal and the binormal vector, the curvature and the torsion of the curve, respectively.\\
Note that, the representation of the arc length assumed parameter $"s"$ was omitted throughout the paper for simplicity, unless otherwise stated.
\section{Tangent Associated Curves} \label{t_assoc}
In this section we will define tangent associated curves such that the tangent vector of a given curve lies on the osculating, normal and rectifying plane of its mate.
Take a differentiable curve, $\alpha$ and denote $\alpha^*$ as its associated mate. Let $\{T^*,N^*,B^*\}$ be the Frenet frame of $\alpha^*$. Then the unit vectors lying on osculating, normal and rectifying plane of the curve $\alpha^*$ are given as
\begin{eqnarray}\label{uniVecs}
\label{oyldz} O^*&=&\frac{aT^*+bN^*}{\sqrt{a^2+b^2}},\\
\label{pyldz} P^*&=&\frac{cN^*+dB^*}{\sqrt{c^2+d^2}},\\
\label{ryldz} R^*&=&\frac{eT^*+fB^*}{\sqrt{e^2+f^2}}
\end{eqnarray}
respectively, where $a,b,c,d,e,f \in R$ are some arbitrary real numbers. 

\theoremstyle{definition}
\begin{definition}
Let $\alpha$ and $\alpha^*$ be two differentiable curves. If the tangent vector, T of $\alpha$ is linearly dependent with the vector, $O^*$, then we name the curve $\alpha^*$ as $T-O^*$ associated curve of $\alpha$.
\end{definition}
\begin{theorem}
Let $\alpha^*$ be $T-O^*$ associated curve of $\alpha$ and $\{T^*,N^*,B^*\}$ be the Frenet frame of $\alpha^*$. The relationship of the corresponding Frenet Frames are given as follows:
\begin{eqnarray}
{T^*}&=& \frac{a}{\sqrt{a^2+b^2}}T+\frac{b}{\sqrt{a^2+b^2}}N \nonumber\\
{N^*} &=& \frac{b}{\sqrt{a^2+b^2}}T-\frac{a}{\sqrt{a^2+b^2}}N \nonumber\\
{B^*}&=& -B \nonumber.
\end{eqnarray}
\end{theorem}
\begin{prof}
Since $\alpha$ and $\alpha^*$ are defined as $T-O^*$ associated curves, we write
\begin{equation}\label{bagla}
\alpha^*=\alpha + \lambda T.
\end{equation}
By differentiating the relation (\ref{bagla}), taking its norm and using the Frenet formulae given in (\ref{freneForms}), we have:
\begin{equation}\label{hamtyldz}
T^*= \frac{(1+\lambda')T+\lambda \kappa N}{\sqrt{(1+\lambda')^2+(\lambda \kappa)^2}}.
\end{equation}
Now taking the second derivative of the equation (\ref{bagla}) and referring again to (\ref{freneForms}) we write
\begin{equation}
{\alpha^*}'' = (\lambda''-\lambda \kappa^2) T + \left((1+\lambda')\kappa + (\lambda \kappa)'\right) N +\lambda \kappa \tau B \nonumber.
\end{equation}
The cross production of ${\alpha^*}'$ and ${\alpha^*}''$ leads us the following form,
\begin{equation}
 {\alpha^*}' \wedge \, \, {\alpha^*}''= \left( \lambda^{2} \kappa^{2}\tau \right) T - \left( ( \lambda' +1) \lambda \kappa  \tau \right) N + \left( ( \lambda'  +1 )  \big(\left( \lambda'  +1 \right) \kappa + (\lambda \kappa)'  \big) -\lambda  \kappa (\lambda'' -\lambda \kappa^{2}) \right) B.
\end{equation}

By calling upon (\ref{defns}), we simply calculate $N^*$, and $B^*$ as
\begin{align}\label{hamnbyldz}
N^* &= -\frac{\lambda  \kappa \left(  \left( \lambda'  +1 \right)\big(  \left( \lambda'  +1 \right) \kappa  + \lambda' \kappa  +\lambda \kappa'   \big) -\lambda  \kappa   \left( \lambda''-\lambda  \kappa^{2}  \right)  \right)  T}{\parallel {\alpha^*}'\parallel \, \, \parallel {\alpha^*}' \wedge \, \, {\alpha^*}''\parallel} \nonumber\\
& \quad +\frac{\left( \lambda'  +1 \right) \bigg(  \left( \lambda' +1 \right) \big(  \left( \lambda'  +1 \right) \kappa  + \lambda' \kappa  +\lambda  \kappa'   \big) -\lambda  \kappa   \left( \lambda''-\lambda  \kappa^{2}   \right)  \bigg) N +\lambda  \kappa  \tau \left( \lambda^{2} \kappa^{2}  + \left( \lambda'  +1 \right) ^{2} \right)  B}{\parallel {\alpha^*}'\parallel \, \, \parallel {\alpha^*}' \wedge \, \, {\alpha^*}''\parallel},\nonumber \\
B^* &= \quad \frac{\lambda^{2} \kappa^{2}\tau T-\lambda  \kappa  \tau(\lambda'+1) N + (\lambda'  +1 )  \big(\left( \lambda'  +1 \right) \kappa + (\lambda \kappa)'  \big) -\lambda  \kappa \left(\lambda'' -\lambda \kappa^{2} \right) B}{\parallel {\alpha^*}' \wedge \, \, {\alpha^*}''\parallel}.
\end{align}

Now, as we defined the curve $\alpha^*$ to be the $T-O^*$ associated curve of $\alpha $ we deduce that \\
$ \qquad <T,T^*>\,=\,<O^*,T^*>$.
By using this deduction and referring both the relation (\ref{oyldz}) and (\ref{hamtyldz}) we write
\[\frac{(1+\lambda')}{\sqrt{(1+\lambda')^2+(\lambda \kappa)^2}}=\frac{a}{\sqrt{a^2+b^2}}.\]
Simple elementary operations on this relation result the following linear ordinary differential equation (ODE), with $b\neq 0$.
\begin{equation}\label{baginti1o}
1+\lambda'=\frac{a}{b}\lambda \kappa.
\end{equation}
When substituted the given ODE, (\ref{baginti1o}) into (\ref{hamtyldz}) we complete the first part of the proof for $T^*$.\\
Similarly, another deductions can be drawn as
$$<T,N^*>\,=\,<O^*,N^*> \text{      and     } <T,B^*>\,=\,<O^*,B^*>=0$$
and using these we write
\begin{align}
\label{baginti2o} &-\frac{\lambda  \kappa \left[  \left( \lambda'  +1 \right)\bigg(  \left( \lambda'  +1 \right) \kappa  + \lambda' \kappa  +\lambda \kappa'   \bigg) -\lambda  \kappa   \left( \lambda''-\lambda  \kappa^{2}  \right)  \right]  }{\parallel {\alpha^*}'\parallel \, \, \parallel {\alpha^*}' \wedge \, \, {\alpha^*}''\parallel}=\frac{b}{\sqrt{a^2+b^2}},&\\
\label{baginti3o} &\lambda^{2} \kappa^{2}\tau =0,&
\end{align}\\
respectively.
Now when substituted the relations (\ref{baginti1o}), (\ref{baginti2o}) and (\ref{baginti3o}) into both (\ref{hamtyldz}) and (\ref{hamnbyldz}) we complete the proof.
\end{prof}
\begin{corollary}
From (\ref{baginti1o}) and (\ref{baginti3o}) $\lambda \neq 0$ that results $\tau=0$. Therefore it can be easily said that the curve $\alpha$ is a planar curve or equivalently there is no a space curve having a T associated partner such that its tangent lies on the osculating plane of its mate.
\end{corollary}

\begin{theorem}
Let $\alpha^*$ be the $T-O^*$ associated curve of $\alpha$. The curvature, $\kappa^*$ and the torsion, $\tau^*$ of $\alpha^*$ are given as follows.

\begin{equation*}
\kappa^*=\frac{b}{\lambda \sqrt{a^2+b^2}}
\end{equation*}
\begin{eqnarray}
\begin{aligned}
\tau^*&= 0 \nonumber
\end{aligned}
\end{eqnarray}
\end{theorem}
\begin{prof}
By using the definitions in (\ref{curvatures}) and the relation (\ref{baginti1o}) with the fact that $\tau=0$ the proof is completed.
\end{prof}

\begin{theorem}\label{solODE}
Let $\alpha^*$ be the $T-O^*$ associated curve of $\alpha$. The distance between the corresponding points of  $\alpha$ and $ \alpha^*$ in $E^3$ is given as follows:
\begin{align}
d({\alpha},\alpha^*)&=\bigg|e^{\int{\frac{a}{b}\kappa}}\bigg[-\int{e^{-\int{\frac{a}{b}\kappa}}}+c_1\bigg]\bigg|
\end{align}
where $c_1$is an integral constant.
\end{theorem}

\begin{prof}
Rewrite (\ref{baginti1o}) as
\begin{equation}
\lambda'-\frac{a}{b}\kappa \lambda=-1.
\end{equation}
By taking $\mu $ as an integrating factor and multiplying the both hand sides of the latter equation by that we get
\begin{equation}\label{multiplied}
 \mu \lambda'-\mu \frac{a}{b}\kappa \lambda =-\mu . 
\end{equation}
From the product rule of the composite form we write\\
\begin{equation}\label{compForm}
(\mu \lambda)' = \mu \lambda'+\mu' \lambda
\end{equation}
and equate the terms of \ref{compForm} with those in the left hand side of the \ref{multiplied}
we find
\[\mu '=-\mu \frac{a}{b} \kappa.\]
The solution for the integrating factor $\mu$ is given with
\[\int{\frac{\mu'}{\mu}} = -\int{ \frac{a}{b} \kappa } \quad \Rightarrow \quad \mu= e^{-\int{\frac{a}{b} \kappa}+c} \]
On the other hand, the use of integrating factor let us to write following relation
\[[\mu \lambda]' = -\mu .\]
Integrating both hand sides of this equation
\[\mu \lambda +c_o =-\int{\mu}  \]
and leaving $\lambda$ all alone we get
\[\lambda=\frac{-\int{\mu}-c_o}{\mu} .\]
By substituting $\mu$ in place, we finally get
\[\lambda= e^{\int{\frac{a}{b}\kappa}}\bigg[-\int{e^{-\int{\frac{a}{b}\kappa}}}+c_1\bigg].\]
\end{prof}

\begin{definition}
Let $\alpha$ and $\alpha^*$ be two differentiable curves. If the tangent vector, T of $\alpha$ is linearly dependent with the vector, $P^*$ defined in (\ref{pyldz}), then we name the curve $\alpha^*$ as $T-P^*$ associated curve of $\alpha$.
\end{definition}

\begin{theorem}
Let $\alpha$ and $\alpha^*$ be $T-P^*$ associated curves. The relationship of the corresponding Frenet Frames are given as follows:

\begin{eqnarray}
T^*&=& N \nonumber\\
N^* &=& \frac{-c}{\sqrt{c^2+d^2}}T+\frac{d}{\sqrt{c^2+d^2}}B \nonumber\\
B^*&=& \frac{d}{\sqrt{c^2+d^2}}T+\frac{c}{\sqrt{c^2+d^2}}B \nonumber.
\end{eqnarray}
\end{theorem}

\begin{prof}
Since we defined the curve $\alpha^*$ to be as $T-P^*$ associated curve of $\alpha $ we could deduce that \\
$<T,N^*>\,=\,<P^*,N^*>$.
Using this, together with the relations (\ref{pyldz}) and (\ref{hamnbyldz}) results the following:
\begin{equation}\label{baginti1p}
-\frac{\lambda  \kappa \left[  \left( \lambda'  +1 \right)\bigg(  \left( \lambda'  +1 \right) \kappa  + \lambda' \kappa  +\lambda \kappa'   \bigg) -\lambda  \kappa   \left( \lambda''-\lambda  \kappa^{2}  \right)  \right]  }{\parallel {\alpha^*}'\parallel \, \,\parallel {\alpha^*}' \wedge \, \, {\alpha^*}''\parallel}=\frac{c}{\sqrt{c^2+d^2}}.
\end{equation}

By the same manner, the following relation can be derived
$<T,B^*>\,=\,<P^*,B^*>$ that reveals
\begin{equation} \label{baginti2p}
 \frac{\lambda^{2} \kappa^{2}\tau}{\parallel {\alpha^*}' \wedge \, \,{\alpha^*}''\parallel}=\frac{d}{\sqrt{c^2+d^2}}.
\end{equation}
Another deduction that $<T,T^*>\,=\,<P^*,T^*>=0 $ provides
\begin{equation}\label{baginti3p}
1+\lambda'=0 \text{ and so } \lambda=-s+c, \qquad \text{ where c is the integral constant}
\end{equation}
Utilizing these three relations, (\ref{baginti1p}), (\ref{baginti2p}) and (\ref{baginti3p}) results what is stated in the theorem.

Note that, by substituting (\ref{baginti3p}) first in both (\ref{baginti1p}) and (\ref{baginti2p}), we find the following relations
\begin{equation}\label{gecisler}
\frac{\kappa}{\sqrt{\kappa^2+\tau^2}}=\frac{-c}{\sqrt{c^2+d^2}} \text{ and } \frac{\tau}{\sqrt{\kappa^2+\tau^2}}=\frac{d}{\sqrt{c^2+d^2}}
\end{equation}
,respectively which points out that $c=-\kappa$ and $d=\tau$ and since by definition $\kappa \geq 0$, $c \leq 0$. 
\end{prof}

\begin{corollary}
It can be easily seen that if $\alpha^*$ is $T-P^*$ associated curve of $\alpha$, then $\alpha^*$ is the involute of $\alpha$. 
\end{corollary}

\begin{theorem}
Let $\alpha^*$ be the $T-P^*$ associated curve of $\alpha$. The curvature, $\kappa^*$ and the torsion, $\tau^*$ of $\alpha^*$ are given as follows.
\begin{align}
\kappa^*&=\frac{\tau \sqrt{c^2+d^2}}{d \lambda \kappa} =\frac{\sqrt{\kappa^2+\tau^2}}{\lambda \kappa} \\
\tau^*&= \frac{\kappa \tau' -\kappa' \tau}{\lambda \kappa (\kappa^2+\tau^2)} \nonumber
\end{align}
\end{theorem}

\begin{prof}
The proof can be done by using \ref{baginti3p} and \ref{gecisler}.
\end{prof}

\begin{theorem} \label{involuteTheo}
Let $\alpha^*$ be the $T-P^*$ associated curve of $\alpha$. The distance between the corresponding points of  $\alpha$ and $ \alpha^*$ in $E^3$ is given as follows:
\begin{equation}
d({\alpha^*},\alpha)=\bigg|-s+c \bigg|
\end{equation}
\end{theorem}

\begin{prof}
The proof is trivial.
\end{prof}

\begin{definition}{}
Let $\alpha$ and $\alpha^*$ be two differentiable curves. If the tangent vector, T of $\alpha$ is linearly dependent with the vector, $R^*$ defined in (\ref{ryldz}), then we name the curve $\alpha^*$ as $T-R^*$ associated curve of $\alpha$.
\end{definition}

\begin{theorem}
Let $\alpha$ and $\alpha^*$ be $T-R^*$ associated curves. The relationship of the corresponding Frenet Frames are given as follows:

\begin{eqnarray}
T^*&=& \frac{e}{\sqrt{e^2+f^2}}T+\frac{f}{\sqrt{e^2+f^2}}N \nonumber\\
N^* &=& B \nonumber\\
B^*&=& \frac{f}{\sqrt{e^2+f^2}}T-\frac{e}{\sqrt{e^2+f^2}}N \nonumber.
\end{eqnarray}
\end{theorem}

\begin{prof}
Since we defined the curve $\alpha^*$ to be as $T-R^*$ associated curve of $\alpha $ we could deduce that \\
$<T,T^*>\,=\,<R^*,T^*>$.
By using this deduction and referring both the relation (\ref{ryldz}) and (\ref{hamtyldz}) we write
\[\frac{(1+\lambda')}{\sqrt{(1+\lambda')^2+(\lambda \kappa)^2}}=\frac{e}{\sqrt{e^2+f^2}},\]
and with some simple elementary operations on this relation we come up with the following linear ordinary differential equation (ODE), with $f\neq 0$.
\begin{equation}
\label{baginti1r} 1+\lambda'=\frac{e}{f}\lambda \kappa.
\end{equation}
When substituted the given ODE into (\ref{hamtyldz}) we complete the first part of the proof for $T^*$.\\
Similarly, another deduction can be drawn as $<T,B^*>\,=\,<R^*,B^*>$  which results
\begin{align}\label{baginti2r}
\frac{\lambda^{2} \kappa^{2}\tau}{\parallel {\alpha^*}' \wedge \, \, {\alpha^*}''\parallel} = \frac{f}{\sqrt{e^2+f^2}}, \text{  and so  } \parallel {\alpha^*}' \wedge \, \, {\alpha^*}''\parallel =  \lambda^{2} \kappa^{2}\tau \frac{\sqrt{e^2+f^2}}{f}.
\end{align}
Now when substituted the relations (\ref{baginti1r}) and (\ref{baginti2r}) into (\ref{hamnbyldz}) we complete the proof for $B^*$.\\
A final inference on the idea of $T-R^*$ association can be drawn as $<T,N^*>\,=\,<R^*,N^*>=0$. This puts the following equation forward
\begin{equation}\label{baginti3r}
-\lambda  \kappa \left[  \left( \lambda'  +1 \right)\bigg(  \left( \lambda'  +1 \right) \kappa  + \lambda' \kappa  +\lambda \kappa'   \bigg) -\lambda  \kappa   \left( \lambda''-\lambda  \kappa^{2}  \right)  \right] =0.
\end{equation}
By substituting (\ref{baginti1r}), (\ref{baginti2r}) and (\ref{baginti3r}) in (\ref{hamnbyldz}) the proof is completed for $N^*$ and all.
\end{prof}

\begin{theorem}
Let $\alpha^*$ be the $T-R^*$ associated curve of $\alpha$ The curvature, $\kappa^*$ and the torsion $\tau^*$ of $\alpha^*$ are given as follows.

\begin{equation}\label{kapayldz}
\kappa^*=\frac{\tau f^2}{\lambda \kappa (e^2+f^2)}
\end{equation}

\begin{equation}\label{tauyldz}
\tau^*= \frac { f\left( \kappa^{2}\tau
  {e}^{3}+ \kappa^{2}\tau  e{f}^{2}+ \tau^{3}e{f}^{2}+\kappa   \tau' {e}^{2}f+\kappa  
\tau' {f}^{3}- \kappa' \tau {e}^{2}f- \kappa' \tau
{f}^{3} \right)}{\lambda  \kappa  \left(  \tau^{2}{e}^{2}{f}^{2}+ \tau^{2}{f}^{4}+ \kappa^{2}{e}^{4}+2\, \kappa^{2}{e}^{2}{f}^{2}+ \kappa^{2}{f}^{4} \right) }.
\end{equation}

\end{theorem}

\begin{prof}
By the definitions given in (\ref{curvatures}) and substituting (\ref{baginti1r}) and (\ref{baginti2r}) into this we easily calculate (\ref{kapayldz}).
On the other hand the third derivative of (\ref{bagla}) is
\begin{equation*}
\begin{split}
{\alpha^*}'''&=\left( -3\, \lambda' \kappa^{2}-3\,\lambda  \kappa  \kappa'  - \kappa^{2}+\lambda'''   \right) T +\left( - \kappa^{3}\lambda  -\lambda  \kappa   \tau^{2}+3\, \lambda'' \kappa  +3\, \lambda' \kappa'  +\lambda  \kappa''  +\kappa'   \right) N \\ \nonumber
&+ \left( 3\, \lambda' \kappa  \tau  +\lambda  \kappa  \tau'  +2\,\lambda \kappa'\tau  +\kappa  \tau   \right) B. \nonumber
\end{split}
\end{equation*}
From (\ref{curvatures}), and using (\ref{baginti1r}) we calculate $\tau$ as given in the theorem.
\end{prof}

\begin{theorem}
Let $\alpha^*$ be the $T-R^*$ associated curve of $\alpha$. The distance between the corresponding points of  $\alpha$ and $ \alpha^*$ in $E^3$ is given as follows:
\begin{equation}
d({\alpha^*},\alpha)=\bigg|e^{\int{\frac{e}{f}\kappa}}\bigg[-\int{e^{-\int{\frac{e}{f}\kappa}}}+c_2\bigg]\bigg|
\end{equation}
\end{theorem}

\begin{prof}
The proof is the same as Theorem (\ref{solODE}).
\end{prof}

\section{Examples}
In this section, we exhibit several examples for tangent associated curves for each three planes.
\begin{enumerate}
\item Let $\alpha$ be chosen a unit speed circle as a planar curve given with a parameterization\\ $\alpha(s)=(cos(s),sin(s),0)$.
Since $\alpha$ is chosen to be a circle $\kappa=1$. By taking $a=b=1$, the general solution for the given ODE in \ref{baginti1o} is
$$\lambda(s) = 1+e^s \cdot c_0$$ where $c_0$ is the integral constant.


\begin{figure}[H]
\centering
\subfigure[$c_0=-1$]{
\includegraphics[width=.375\textwidth]{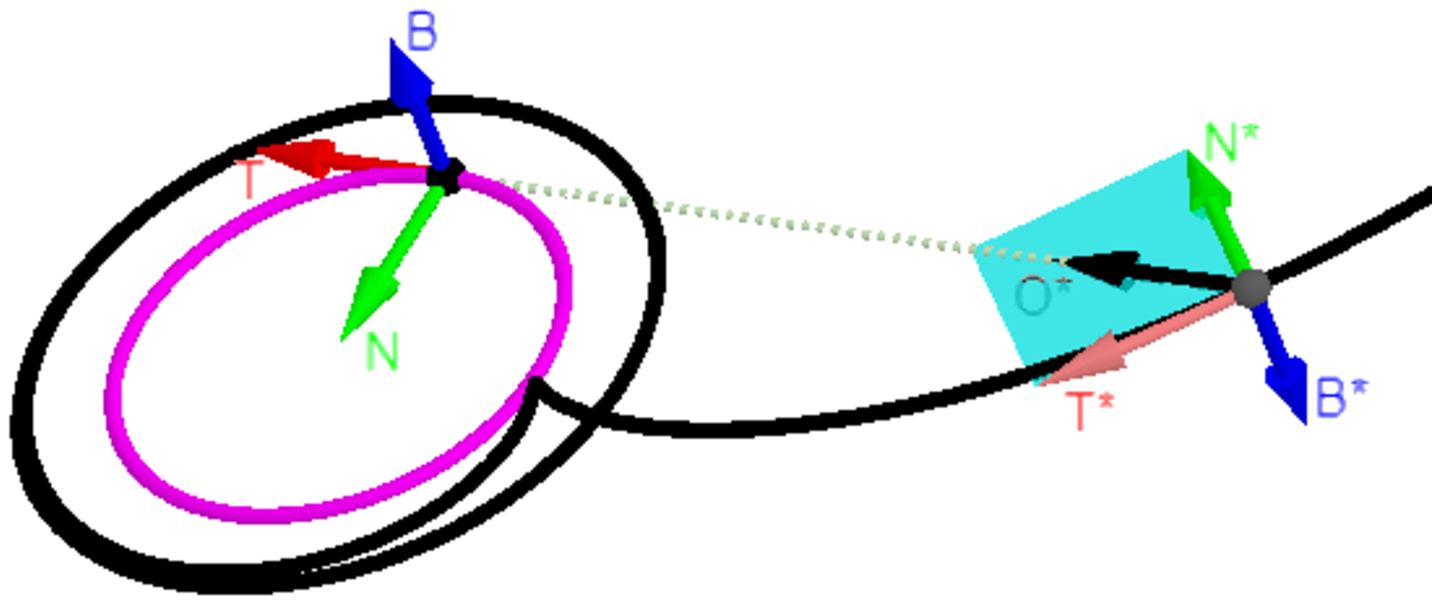}
}
\subfigure[$c_0=0$]{
\includegraphics[width=.225\textwidth]{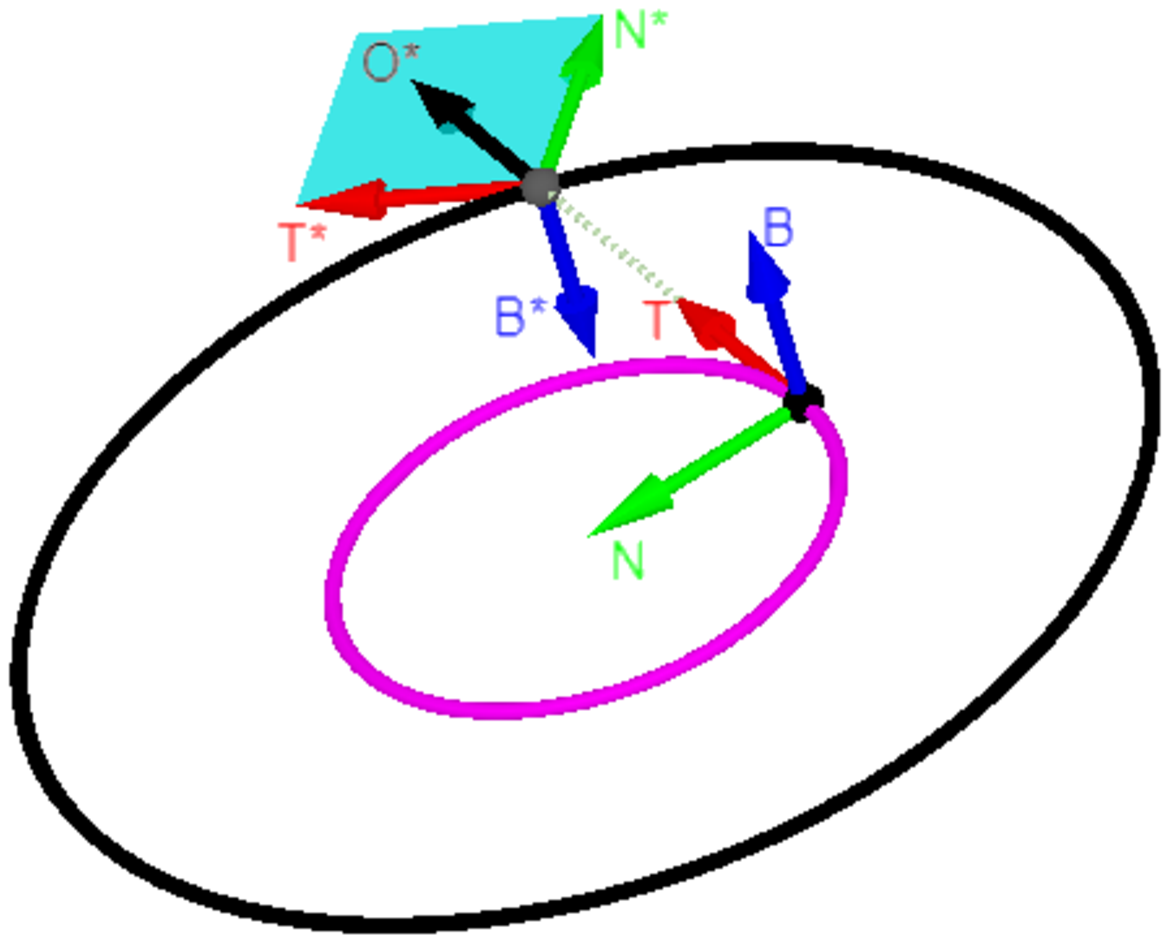}
}
\subfigure[$c_0=1$]{
\includegraphics[width=.250\textwidth]{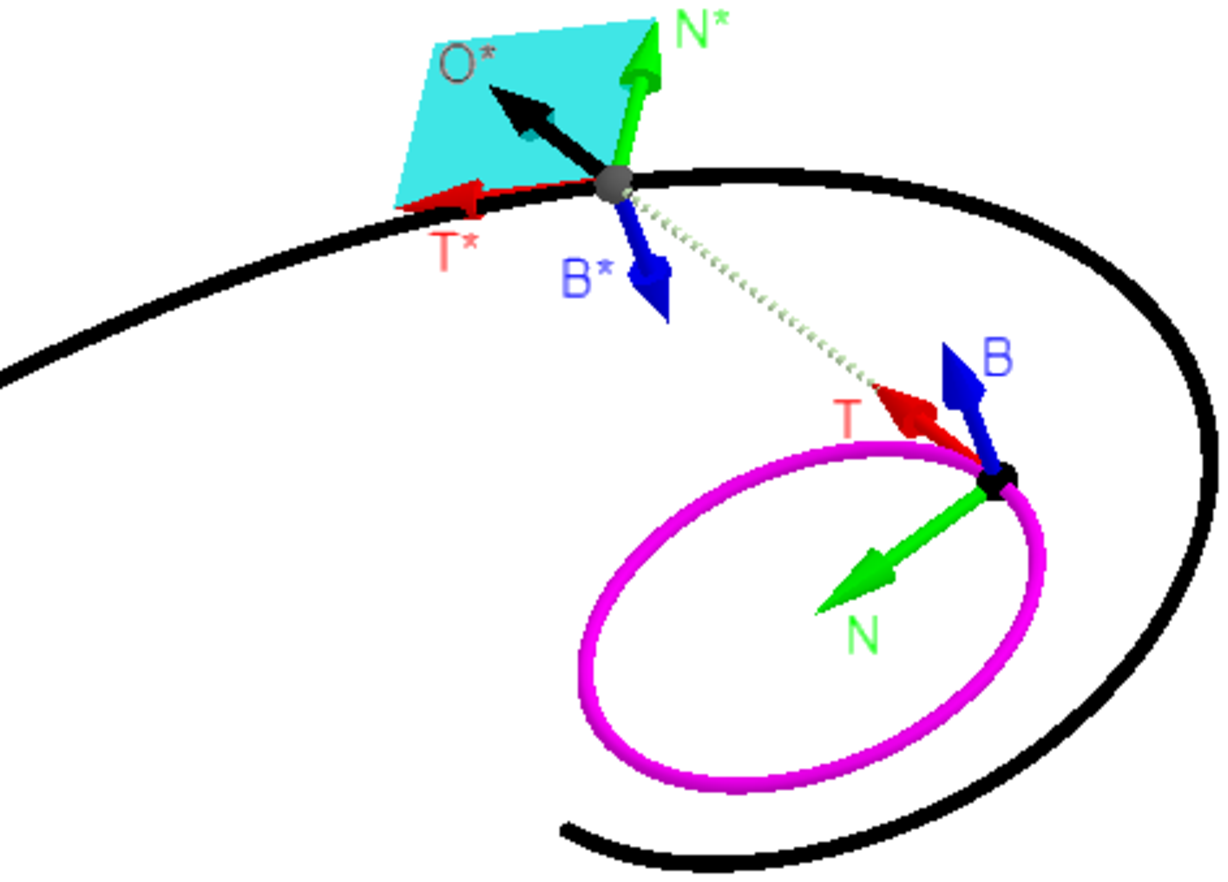}
}
\caption{The curve $\alpha$ (pink) and its $T-O^*$ associated mate $\alpha^*$ (black)}
\end{figure}

\item Let $\alpha$ be chosen a unit speed helix given with a parameterization $\alpha(s)=\frac{1}{\sqrt{2}}(cos(s),sin(s),s)$.
Since $\kappa=\tau=\frac{1}{\sqrt{2}}$, the vector $P^*$ should be formed by the values of $c$ and $d$ such that $-c=d=\frac{1}{\sqrt{2}}$. From theorem (\ref{involuteTheo}) we write
$$\lambda(s) = -s + c_0$$ where $c_0$ is the integral constant.


\begin{figure}[H]
\centering
\subfigure[$c_0=-1$]{
\includegraphics[width=.275\textwidth]{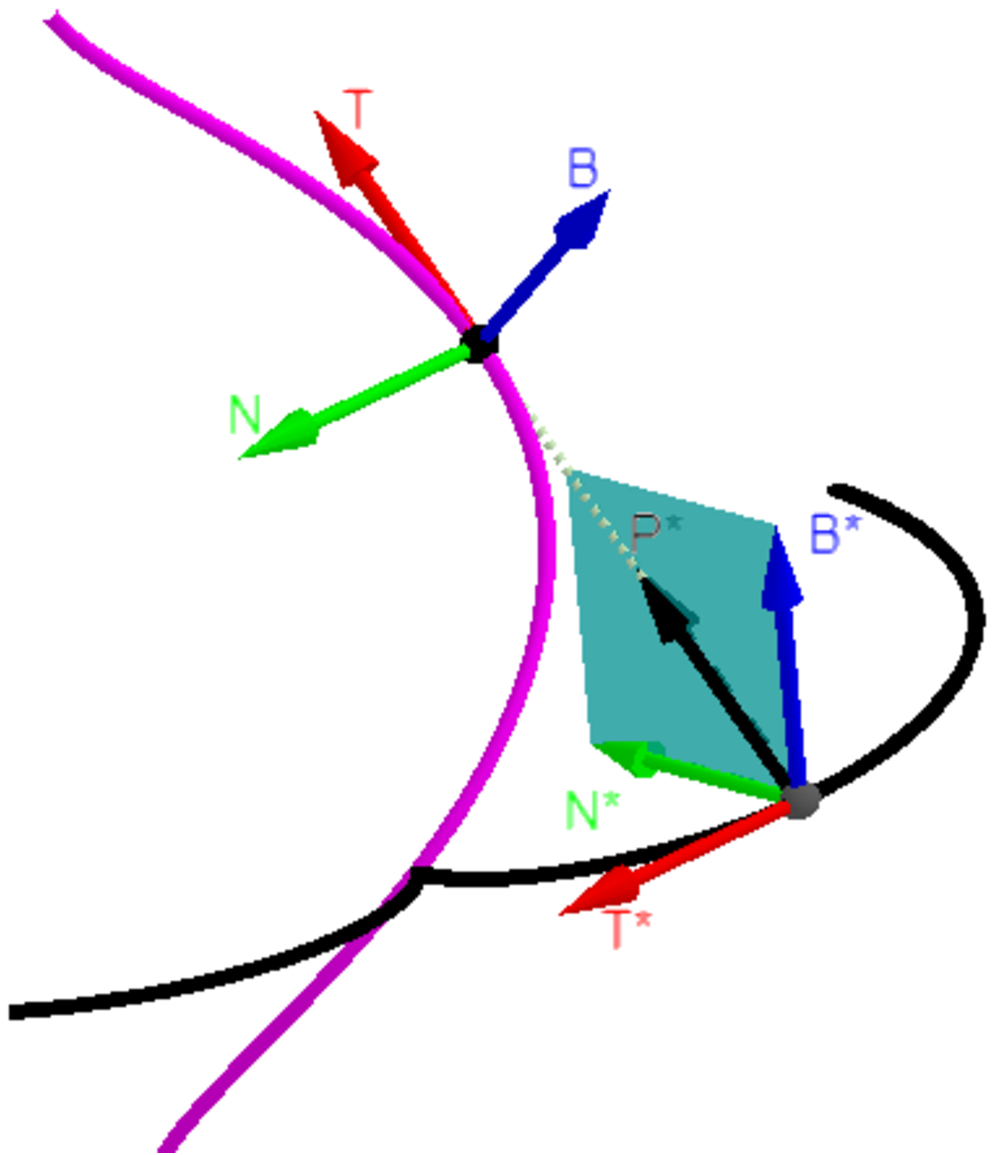}
}
\subfigure[$c_0=0$]{
\includegraphics[width=.275\textwidth]{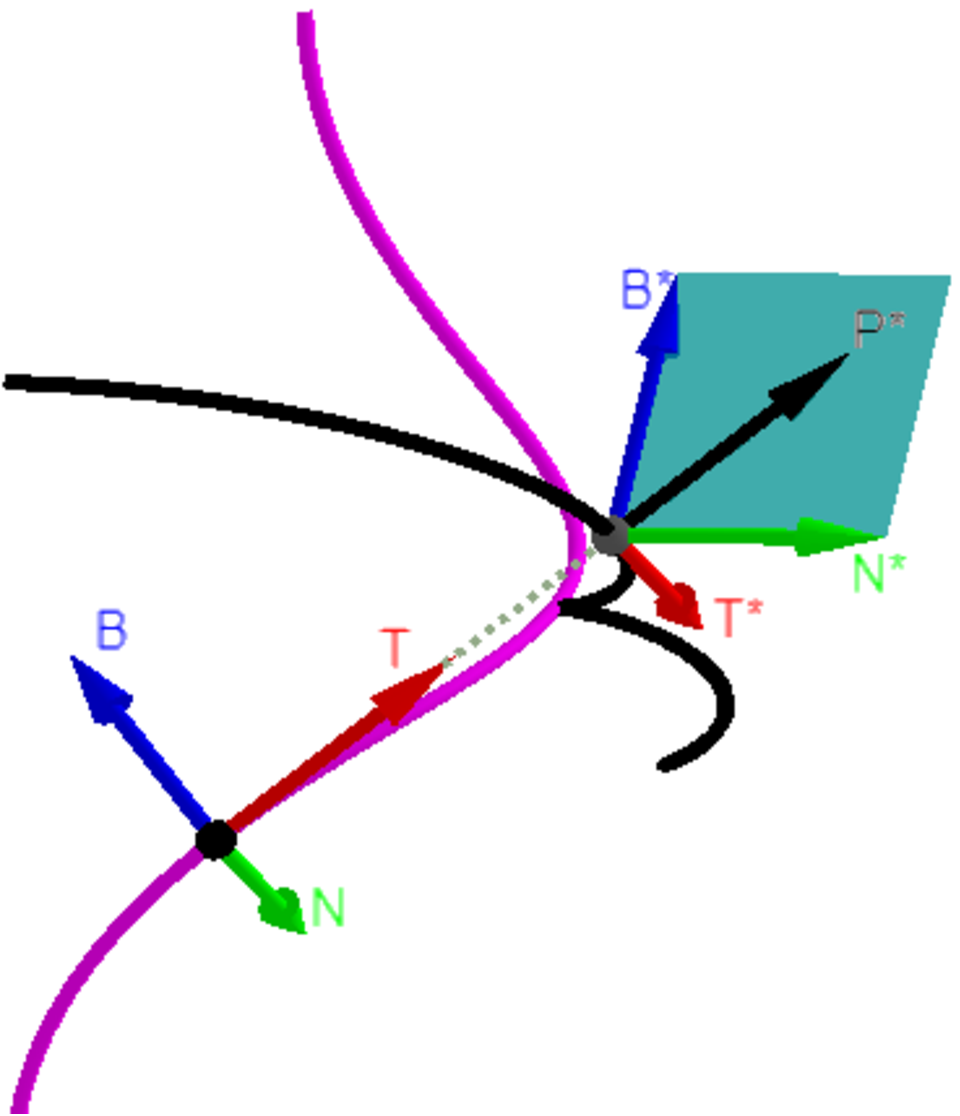}
}
\subfigure[$c_0=1$]{
\includegraphics[width=.275\textwidth]{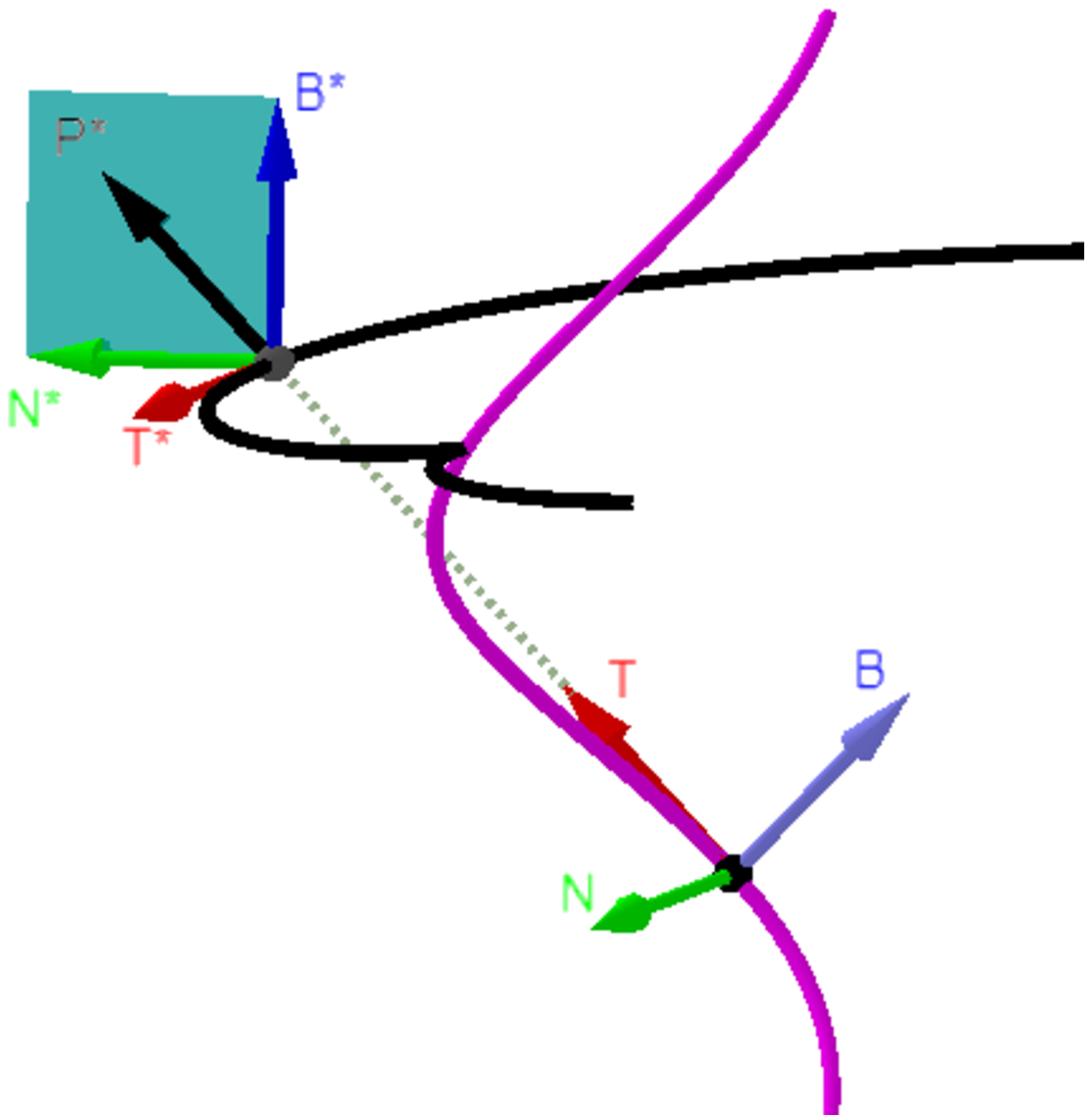}
}
\caption{The curve $\alpha$ (pink) and its $T-P^*$ associated mate $\alpha^*$ (black)}
\end{figure}

\item By referring the same curve given in (ii) we know that $\kappa=\frac{1}{\sqrt{2}}$. The general solution for the ODE in \ref{baginti1o} is this time 
$$\lambda(s) = \sqrt{2}+e^{\frac{\sqrt{2}}{2}s} \cdot c_0$$ where $c_0$ is the integral constant.
for $e=f=1$

\begin{figure}[H]
\centering
\subfigure[$c_0=-1$]{
\includegraphics[width=.375\textwidth]{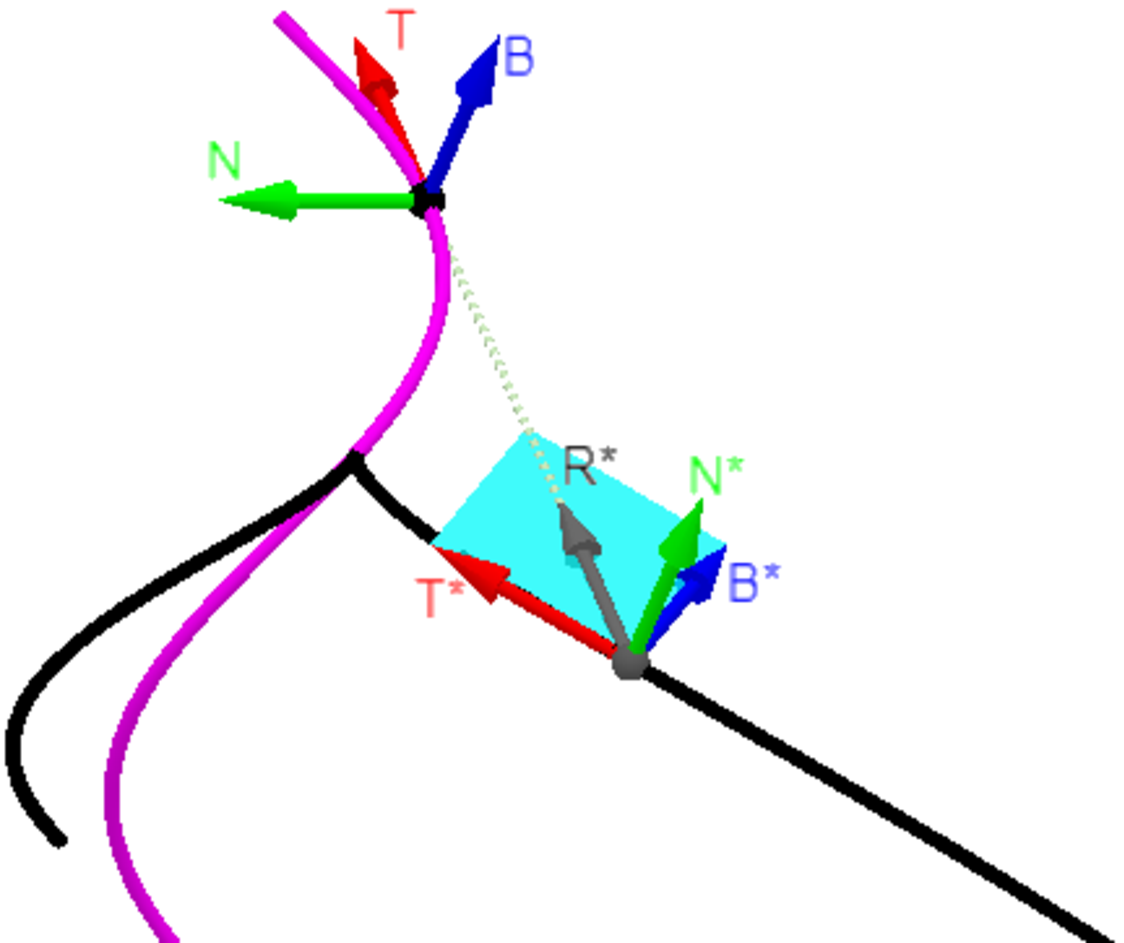}
}
\subfigure[$c_0=0$]{
\includegraphics[width=.225\textwidth]{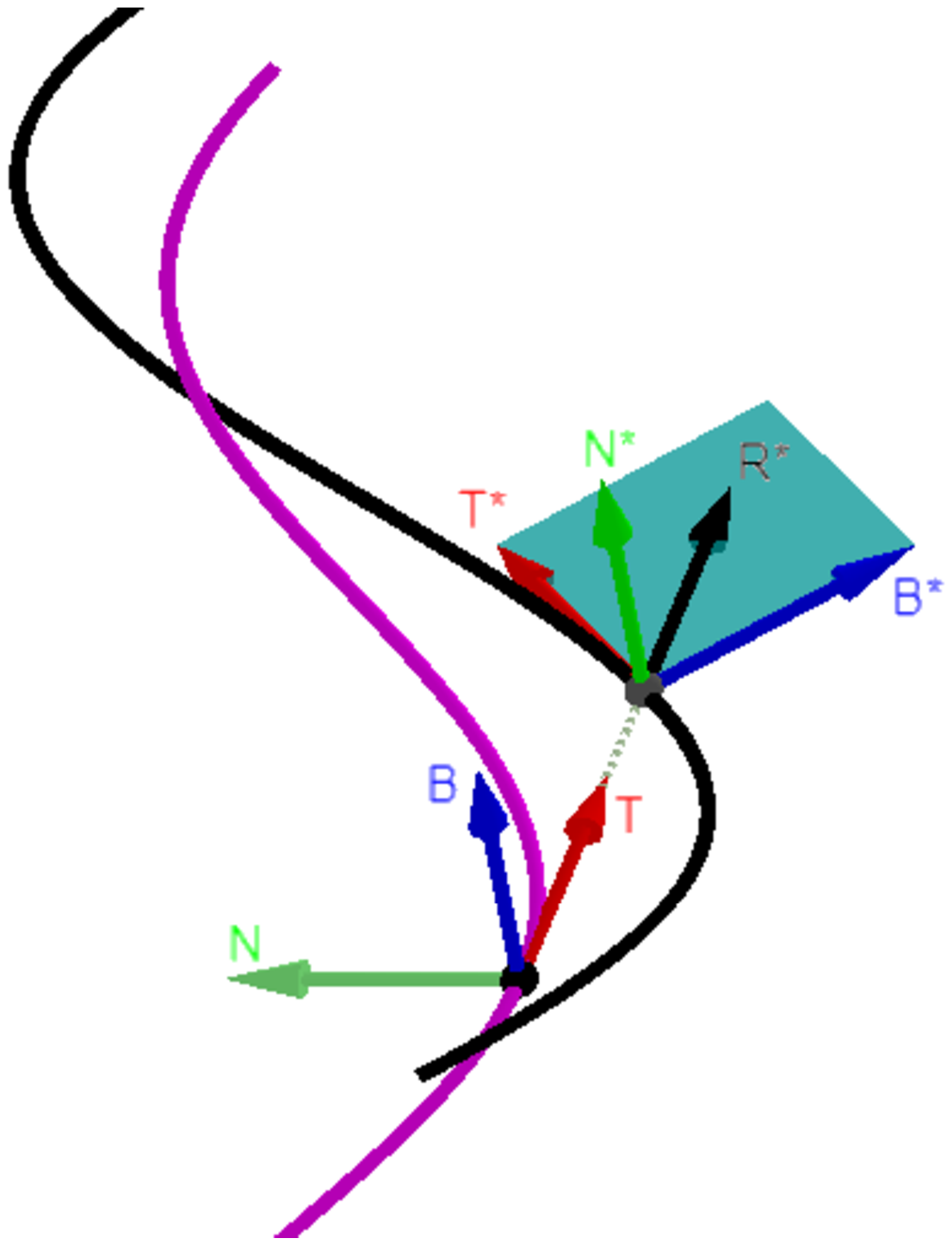}
}
\subfigure[$c_0=1$]{
\includegraphics[width=.250\textwidth]{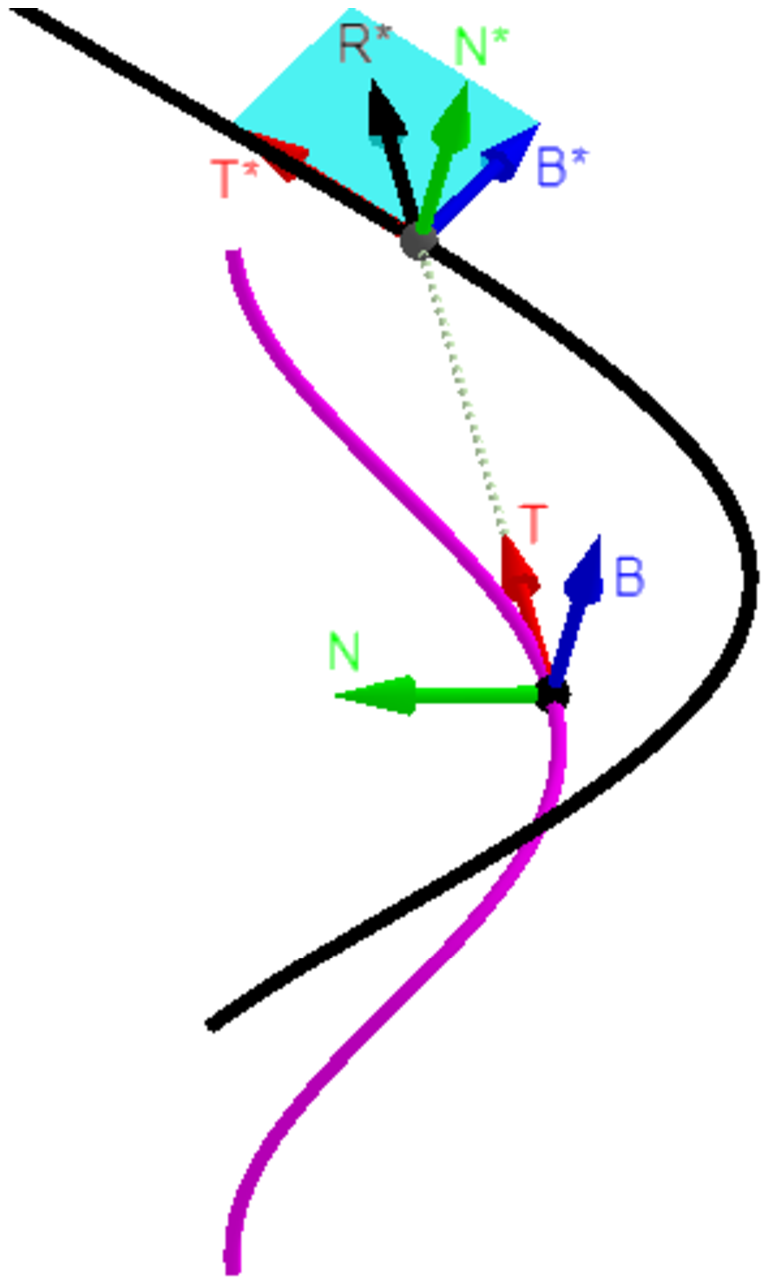}
}
\caption{The curve $\alpha$ (pink) and its $T-R^*$ associated mate $\alpha^*$ (black)}
\end{figure}

\end{enumerate}

\section{Principal Normal Associated Curves}
In this section, we define principal normal associated curves such that the principal normal vector of a given curve lies on the osculating, normal and rectifying plane of its mate.
\theoremstyle{definition}
\begin{definition}
Let $\alpha$ and $\alpha^*$ be two differentiable curves. If the principal normal, N of $\alpha$ is linearly dependent with the vector, $O^*$, then we name the curve $\alpha^*$ as $N-O^*$ associated curve of $\alpha$.
\end{definition}

\begin{theorem}
Let $\alpha^*$ be $N-O^*$ associated curve of $\alpha$ and $\{T^*,N^*,B^*\}$ be the Frenet frame of $\alpha^*$. The relationship of the corresponding Frenet Frames are given as follows:
\begin{align}
{T^*}&=\frac{1}{\sqrt{a^2+b^2}} \left( {\frac { \left(- \lambda \kappa  +1 \right) b}{\sqrt {(- \lambda \kappa  +1 )^2 +(\lambda \tau)^2 }}}T +a N +{\frac {\lambda  \tau b}{\sqrt {(- \lambda \kappa  +1 )^2 +(\lambda \tau)^2 }}}B \right) \nonumber\\
{N^*} &=\frac{b}{\sqrt{a^2+b^2}}\left( \frac{- \mathbf{M} }{\mathbf{M}  \left(\lambda  \kappa  -1 \right)-\mathbf{K}  \lambda  \tau}T+ N + \frac{\mathbf{K} }{\mathbf{M}  \left(\lambda  \kappa  -1 \right)-\mathbf{K}  \lambda  \tau}B \right) \nonumber\\
{B^*}&=\frac{b(\mathbf{K}  T+ \mathbf{M} B)}{a(\mathbf{M}  \left(\lambda  \kappa  -1 \right)-\mathbf{K}  \lambda  \tau)} \nonumber,
\end{align}
where the coefficients $\mathbf{K} \text{ and } \mathbf{M}$ are
\begin{align}
\mathbf{K}&=\lambda' \left( \lambda  \tau'  +2\, \lambda' \tau   \right) -\lambda  \tau   \left(  \left( -\lambda  \kappa  +1 \right) \kappa  -\lambda \tau^{2}+\lambda''   \right)\nonumber \\
\mathbf{M}&=  \left( -\lambda  \kappa  +1 \right)  \left(  \left( -\lambda  \kappa  +1 \right) \kappa  -\lambda \tau^{2}+\lambda''   \right) - \lambda'  \left( -\lambda \kappa'  -2\, \lambda'\kappa   \right)\nonumber.
\end{align}
\end{theorem}

\begin{prof}
Since $\alpha$ and $\alpha^*$ are defined as $N-O^*$ associated curves, we write
\begin{equation}\label{bagla2}
\alpha^*=\alpha + \lambda N.
\end{equation}
By differentiating the relation (\ref{bagla2}), using the Frenet formulae given in (\ref{freneForms}) and taking the norm, we have:
\begin{equation}\label{hamtyldz2}
T^*= \frac {(-\lambda \kappa  +1)T+\lambda' N+ \lambda  \tau B}{\sqrt {(-\lambda \kappa  +1)^2+(\lambda' )^2+(\lambda  \tau)^2}}
\end{equation}
Next taking the second derivative of the equation (\ref{bagla2}) and referring again to (\ref{freneForms}) result the following relation.
\begin{equation}
{\alpha^*}'' =  \left( -\lambda \kappa'  -2\, \lambda' \kappa   \right)T  + \left(  \left( -\lambda  \kappa  +1 \right) \kappa  -\lambda \tau^{2}+\lambda''   \right) N  + \left( \lambda \tau'  +2\, \lambda' \tau   \right) B \nonumber.
\end{equation}

The cross production of ${\alpha^*}'$ and ${\alpha^*}''$ leads us the following form,
\begin{equation}
{\alpha^*}' \wedge \, \, {\alpha^*}''= \mathbf{K}  T+ \mathbf{L} N + \mathbf{M} B \nonumber
\end{equation}
where  $\mathbf{K}, \quad \mathbf{L} \quad \text{and} \quad \mathbf{M}$ are assigned to be as
\begin{align}\label{kisaltmalar1}
\mathbf{K}&=\lambda' \left( \lambda  \tau'  +2\, \lambda' \tau   \right) -\lambda  \tau   \left(  \left( -\lambda  \kappa  +1 \right) \kappa  -\lambda \tau^{2}+\lambda''   \right)\nonumber \\
\mathbf{L}&= \left( -\lambda  \kappa  +1 \right)  \left( \lambda  \tau' +2\, \lambda' \tau   \right) +\lambda  \tau   \left( -\lambda  \kappa'  -2\, \lambda' \kappa   \right)\\
\mathbf{M}&=  \left( -\lambda  \kappa  +1 \right)  \left(  \left( -\lambda  \kappa  +1 \right) \kappa  -\lambda \tau^{2}+\lambda''   \right) - \lambda'  \left( -\lambda \kappa'  -2\, \lambda'\kappa   \right)\nonumber
\end{align}
for the sake of simplicity. Note that the norm, $\parallel {\alpha^*}' \wedge \, \, {\alpha^*}''\parallel=\sqrt{\mathbf{K}^2+\mathbf{L}^2+\mathbf{M}^2}$.\\ 
By referring again the definitions given by (\ref{defns}), we simply calculate $N^*$, and $B^*$ as
\begin{align}\label{hamnbyldz2}
N^* &= \frac{ (\mathbf{L} \lambda \tau  - \mathbf{M} \lambda' )T +(\mathbf{M}  \left(\lambda  \kappa  -1 \right)-\mathbf{K}  \lambda  \tau )N +(\mathbf{K} \lambda'  -\mathbf{L}  \left( -\lambda  \kappa  +1 \right))B}{\sqrt {(-\lambda \kappa  +1)^2+(\lambda' )^2+(\lambda  \tau)^2} \sqrt{\mathbf{K}^2+\mathbf{L}^2+\mathbf{M}^2}}\\
B^* &= \frac{\mathbf{K}  T+ \mathbf{L} N + \mathbf{M} B}{\sqrt{\mathbf{K}^2+\mathbf{L}^2+\mathbf{M}^2}}\nonumber
\end{align}
The intuitive idea is as same as before. Since we defined $\alpha^*$ to be as the $N-O^*$ associated curve of $\alpha $ we can write that $ <N,T^*>\,=\,<O^*,T^*>$.
By using this together with the relations (\ref{oyldz}) and (\ref{hamtyldz2}) we write
\begin{equation}\label{baginti1no}
\frac{\lambda' }{\sqrt {(-\lambda \kappa  +1)^2+(\lambda' )^2+(\lambda  \tau)^2} }=\frac{a}{\sqrt{a^2+b^2}}.
\end{equation}
Similarly, we can write $<N,N^*>\,=\,<O^*,N^*> $ which results the following
\begin{equation}\label{baginti2no} 
\frac{\mathbf{M}  \left(\lambda  \kappa  -1 \right)-\mathbf{K}  \lambda  \tau}{\sqrt {(-\lambda \kappa  +1)^2+(\lambda' )^2+(\lambda  \tau)^2} \sqrt{\mathbf{K}^2+\mathbf{L}^2+\mathbf{M}^2}}=\frac{b}{\sqrt{a^2+b^2}}.
\end{equation}
and by the same idea $<N,B^*>\,=\,<O^*,B^*>=0$ we get
\begin{equation}\label{baginti3no} 
\frac{\mathbf{L}}{\sqrt{\mathbf{K}^2+\mathbf{L}^2+\mathbf{M}^2}}=0.
\end{equation}
When substituted the given three relations (\ref{baginti1no}), (\ref{baginti2no}) and (\ref{baginti3no}) into (\ref{hamtyldz2}) and (\ref{hamnbyldz2}), we complete the proof.
\end{prof}
Note that none of the differential equations given above is solvable analytically. However we might solve them under some assumptions. 
\begin{corollary}\label{sonuc1}
If the curve $\alpha$ is helix, then by the relation (\ref{baginti1no}) we get
\[ \lambda'\left( \lambda'' -\frac{a^2}{b^2}  \bigg((\lambda  \kappa  -1) \kappa +\lambda \tau^2\bigg) \right)=0,\]
which is solvable analytically in two folds. First, $\lambda'=0$ corresponding to that $\lambda$ is a constant. This is true iff $a=0$, and if $a=0$ then $O^*=N^*$ that is clearly the definition of Bertrand curves.\\
When considered the second factor of the latter relation we come up with a non homogeneous linear second order differential equation with constant coefficients. There we have a complex solution for this case which is given as
\[\lambda =\sin \left( {\frac {a i \sqrt {{\kappa}^{2}+{
\tau}^{2}}}{b}} \right) {c_1}+\cos \left( {\frac {a i \sqrt {{
\kappa}^{2}+{\tau}^{2}}}{b}} \right) {c_2}+{\frac {\kappa}{{
\kappa}^{2}+{\tau}^{2}}}, \qquad \quad i^2=-1.\]
Since $sin(i x)=isinh(x) \text{  and  } cos(i x)=cosh(x)$, we reform the given relation as 
\[\lambda =i \sinh \left( {\frac {a \sqrt {{\kappa}^{2}+{
\tau}^{2}}}{b}} \right) {c_1}+\cosh \left( {\frac {a \sqrt {{
\kappa}^{2}+{\tau}^{2}}}{b}} \right) {c_2}+{\frac {\kappa}{{
\kappa}^{2}+{\tau}^{2}}}\]
where $c_1$ and $c_2$ are integration constants.\\
Now, by recalling the relations (\ref{kisaltmalar1}) and (\ref{baginti3no}) under the assumption that $\alpha$ is helix, then we have
\[ \lambda' \tau (-2\lambda \kappa +1)=0.\]
%
This results that $\lambda$ is a constant of the form, $\displaystyle \lambda=\frac{1}{2\kappa}.$
\end{corollary}

\begin{theorem}
Let $\alpha^*$ be the $N-O^*$ associated curve of $\alpha$. The curvature, $\kappa^*$ and the torsion, $\tau^*$ of $\alpha^*$ are given as follows.
\begin{equation}
\begin{aligned}
\kappa^*&=\frac{a^4 (\mathbf{M}  \left(\lambda  \kappa  -1 \right)-\mathbf{K}  \lambda  \tau)}{b(a^2+b^2) (\lambda')^4}\\
\tau^*&= \frac{b \lambda'}{a(\mathbf{M}  \left(\lambda  \kappa  -1 \right)-\mathbf{K}  \lambda  \tau)} \bigg( \mathbf{K} \left( \lambda {\kappa}^{3}+\lambda \kappa {\tau}^{2}-3 \lambda' \kappa' -\lambda \kappa'' -3 \lambda'' \kappa-{\kappa}^{2} \right) \\
&\qquad \qquad \qquad \qquad \qquad \qquad +\mathbf{M  } \left(\kappa \tau -\lambda \kappa^{2} \tau-\lambda \tau^{3} +3  \lambda' \tau' +\lambda \tau'' +3 \lambda'' \tau\right)  \bigg). \nonumber
\end{aligned}
\end{equation}
\end{theorem}

\begin{prof}
By taking the third derivative of (\ref{bagla2}) and using Frenet formulas, we have 
\begin{align}\label{turev3n}
{\alpha^*}'''&=\left( \lambda {\kappa}^{3}+\lambda \kappa {\tau}^{2}-3 \lambda' \kappa' -\lambda \kappa'' -3 \lambda'' \kappa-{\kappa}^{2} \right) T\nonumber \\ 
&+\left( \lambda''' -3 \lambda' ({\kappa}^{2}+{\tau}^{2})-3\lambda( \kappa \kappa'+ \tau \tau' )+\kappa' \right) N \\ 
&+ \left(\kappa \tau -\lambda \kappa^{2} \tau-\lambda \tau^{3} +3  \lambda' \tau' +\lambda \tau'' +3 \lambda'' \tau\right) B. \nonumber
\end{align}
Now we recall the relations (\ref{baginti1no}), (\ref{baginti2no}) and (\ref{baginti3no}), and substitute these into the definitions given in (\ref{curvatures}), we complete the proof.
\end{prof}

\theoremstyle{definition}
\begin{definition}
Let $\alpha$ and $\alpha^*$ be two differentiable curves. If the principal normal, N of $\alpha$ is linearly dependent with the vector, $P^*$, then we name the curve $\alpha^*$ as $N-P^*$ associated curve of $\alpha$.
\end{definition}

\begin{theorem}
Let $\alpha^*$ be $N-P^*$ associated curve of $\alpha$ and $\{T^*,N^*,B^*\}$ be the Frenet frame of $\alpha^*$. The relationship of the corresponding Frenet Frames are given as follows:
\begin{align}
{T^*}&= {\frac {-\lambda\,\kappa  +1}{\sqrt {(- \lambda \kappa  +1 )^2 +(\lambda \tau)^2 }}}T+ {\frac {\lambda\,\tau  }{\sqrt {(- \lambda \kappa  +1 )^2 +(\lambda \tau)^2 }}}B \nonumber\\
{N^*} &=\frac{1}{\sqrt{c^2+d^2}}\left( \frac{ d\lambda \tau}{\sqrt {(- \lambda \kappa  +1 )^2 +(\lambda \tau)^2 }}T+ c N + \frac{-d(-\lambda \kappa +1)}{\sqrt {(- \lambda \kappa  +1 )^2 +(\lambda \tau)^2 }}B \right) \nonumber\\
{B^*} &=\frac{d}{\sqrt{c^2+d^2}}\left( -\frac {\tau   \left( \lambda\,\tau^{2}+ \kappa^{2}\lambda-\kappa   \right) }{\lambda\,\tau \kappa'  -\kappa \tau'\lambda+\tau'  }T+ N + \frac { \left( \lambda\,\kappa  -1 \right)  \left( \lambda\, \tau^{2}+ \left( \kappa   \right) ^{2}\lambda-\kappa   \right) }{\lambda\, \left( -\lambda\,\tau  \kappa'  +\kappa  \tau' \lambda-\tau' \right) } B \right) \nonumber
\end{align}
\end{theorem}

\begin{prof}
Now, since again we defined $\alpha^*$ to be as the $N-P^*$ associated curve of $\alpha $ we can write three of our associative relations as usual which are
\begin{itemize}
\item $ <N,N^*>\,=\,<P^*,N^*>$,
\item $ <N,B^*>\,=\,<P^*,B^*>$,
\item $ <N,T^*>\,=\,<P^*,T^*> =0$.
\end{itemize}
These relations this time result the following three equations
\begin{align}\label{baginti1np}
\bullet \qquad&\frac{\mathbf{M}  \left(\lambda  \kappa  -1 \right)-\mathbf{K}  \lambda  \tau }{\sqrt {(-\lambda \kappa  +1)^2+(\lambda' )^2+(\lambda  \tau)^2} \sqrt{\mathbf{K}^2+\mathbf{L}^2+\mathbf{M}^2}}=\frac{c}{\sqrt{c^2+d^2}},\nonumber \\
\bullet \qquad&\frac{\mathbf{L}}{\sqrt{\mathbf{K}^2+\mathbf{L}^2 +\mathbf{M}^2}}=\frac{d}{\sqrt{c^2+d^2}},\\
\bullet \qquad&\frac{\lambda' }{\sqrt {(-\lambda \kappa  +1)^2+(\lambda' )^2+(\lambda  \tau)^2} }=0. \nonumber
\end{align}

When substituted the latter relations in (\ref{hamtyldz2}) and (\ref{hamnbyldz2}) we complete the proof.
\end{prof}

\begin{corollary}
Note that the third relation in (\ref{baginti1np}) results that $\lambda$ is constant. What we know from literature is that a Bertrand curve has a constant distance as well as the Mannheim curves (see \cite{bertrand}, \cite{liuWang} and \cite{mannheim}). For Bertrand curves we also know that curves share the principal normal vectors as common, on the other hand for Mannheim curves, they share the property of the parallelization of principal normal and binormal vectors. By our result, we see that if the principal normal vector of any given curve coincides the unit vector spanned by principal normal and binormal vectors of its mate, then the distance of two curves is constant, in general.
\end{corollary}
\begin{theorem}
Let $\alpha^*$ be the $N-P^*$ associated curve of $\alpha$. The curvature, $\kappa^*$ and the torsion, $\tau^*$ of $\alpha^*$ are given as follows.
\begin{equation}
\begin{aligned}
\kappa^*&=\frac{\mathbf{l}\,c^3 \sqrt{c^2+d^2}}{d^4 \left(\mathbf{m}  (\lambda  \kappa -1 )-\mathbf{k} \lambda  \tau \right)^3}\\
\tau^*&= \frac{\mathbf{l}^2 (c^2+d^2)}{d^2} \bigg( \mathbf{k} \left(  \kappa^{3}\lambda+\kappa \tau^{2}\lambda-\lambda\,\kappa'' - \kappa^{2} \right)+\mathbf{l} \left( -3\,\lambda\,\tau  \tau'  -3\,\lambda\, \kappa' \kappa  +\kappa'   \right) \\
&\qquad \qquad \qquad +\mathbf{m} \left( -\kappa^{2}\tau  \lambda- \tau^{3}\lambda+\lambda\,\tau''  +\kappa  \tau   \right)\bigg) \nonumber
\end{aligned}
\end{equation}
where $\mathbf{k},\,\mathbf{l},\text{ and  }\mathbf{m}$ are the coefficients of which $\mathbf{K},\,\mathbf{L},\text{ and  }\mathbf{M}$ reformed with $\lambda'=0$, respectively.
\end{theorem}

\begin{prof}
By referring the relations (\ref{baginti1np}) together with (\ref{turev3n}), the proof is trivial. 
\end{prof}

\theoremstyle{definition}
\begin{definition}
Let $\alpha$ and $\alpha^*$ be two differentiable curves. If the principal normal, N of $\alpha$ is linearly dependent with the vector, $R^*$, then we name the curve $\alpha^*$ as $N-R^*$ associated curve of $\alpha$.
\end{definition}

\begin{theorem}
Let $\alpha^*$ be $N-R^*$ associated curve of $\alpha$ and $\{T^*,N^*,B^*\}$ be the Frenet frame of $\alpha^*$. The relationship of the corresponding Frenet Frames are given as follows:
\begin{equation}
\begin{aligned}
{T^*}&=\frac{1}{\sqrt{e^2+f^2}} \left( {\frac { \left(- \lambda \kappa  +1 \right) f}{\sqrt {(- \lambda \kappa  +1 )^2 +(\lambda \tau)^2 }}}T +e N +{\frac {\lambda  \tau f}{\sqrt {(- \lambda \kappa  +1 )^2 +(\lambda \tau)^2 }}}B \right), \nonumber \\
{N^*} &=\frac{ef}{\sqrt{e^2+f^2}}\bigg( \bigg( \frac{\lambda \tau}{\lambda'} -{\frac { \left( -\lambda  \kappa  +1 \right)  \left(  \left( -\lambda  \kappa  +1 \right) \kappa  -\lambda \tau^{2}+\lambda''   \right) - \lambda' \left( -\lambda \kappa'  -2\,\lambda' \kappa   \right) }{\left(\lambda  \kappa  -1 \right)  \left( \lambda \tau'  +2\, \lambda'\tau   \right) +\lambda  \tau   \left( -\lambda \kappa'  -2\,\lambda'\kappa   \right) }} \bigg)T \nonumber \\
& \qquad \qquad \qquad + \bigg( {\frac { \lambda' \left( \lambda  \tau' +2\, \lambda' \tau  \right) -\lambda  \tau   \left(  \left( -\lambda  \kappa  +1 \right) \kappa  -\lambda \tau^{2}+\lambda''   \right) }{\left( \lambda  \kappa  -1 \right)  \left( \lambda \tau'  +2\, \lambda' \tau   \right) +\lambda  \tau   \left( -\lambda \kappa'  -2\,\lambda' \kappa   \right) } + \frac{\lambda \kappa-1}{\lambda'}}\bigg)B \bigg), \nonumber \\
{B^*} &=\frac{f}{\sqrt{e^2+f^2}}\bigg( \bigg( \frac { \lambda' \left( \lambda  \tau' +2\, \lambda' \tau  \right) -\lambda  \tau   \left(  \left( -\lambda  \kappa  +1 \right) \kappa  -\lambda \tau^{2}+\lambda''   \right) }{\left( \lambda  \kappa  -1 \right)  \left( \lambda \tau'  +2\, \lambda' \tau   \right) +\lambda  \tau   \left( -\lambda \kappa'  -2\,\lambda' \kappa   \right) } \bigg) T + N \\
& \qquad \qquad  \qquad + \bigg( {\frac { \left( -\lambda  \kappa  +1 \right)  \left(  \left( -\lambda  \kappa  +1 \right) \kappa  -\lambda \tau^{2}+\lambda''   \right) - \lambda' \left( -\lambda \kappa'  -2\,\lambda' \kappa   \right) }{\left(\lambda  \kappa  -1 \right)  \left( \lambda \tau'  +2\, \lambda'\tau   \right) +\lambda  \tau   \left( -\lambda \kappa'  -2\,\lambda'\kappa   \right) }}  \bigg)B   \bigg) 
\end{aligned}
\end{equation}

\end{theorem}
\begin{prof}
Now, since again we defined $\alpha^*$ to be as the $N-R^*$ associated curve of $\alpha $ we can write three of our associative relations as usual which are
\begin{itemize}
\item $ <N,T^*>\,=\,<R^*,T^*>$,
\item $ <N,B^*>\,=\,<R^*,B^*>$,
\item $ <N,N^*>\,=\,<R^*,N^*> =0$.
\end{itemize}
By using these we get
 \begin{align}\label{baginti1nr}
\bullet \qquad &\frac{\lambda' }{\sqrt {(-\lambda \kappa  +1)^2+(\lambda' )^2+(\lambda  \tau)^2} }=\frac{e}{\sqrt{e^2+f^2}}, \nonumber \\
\bullet \qquad &\frac{\mathbf{L}}{\sqrt{\mathbf{K}^2+\mathbf{L}^2+\mathbf{M}^2}}=\frac{f}{\sqrt{e^2+f^2}}, \\
\bullet \qquad &\frac{\mathbf{M}  \left(\lambda  \kappa  -1 \right)-\mathbf{K}  \lambda  \tau}{\sqrt {(-\lambda \kappa  +1)^2+(\lambda' )^2+(\lambda  \tau)^2} \sqrt{\mathbf{K}^2+\mathbf{L}^2+\mathbf{M}^2}}=0.\nonumber
\end{align}
When substituted the above expressions into (\ref{hamtyldz2}) and (\ref{hamnbyldz2}) the proof is complete.
\end{prof}
\begin{corollary}
The only analytically solvable equation in (\ref{baginti1nr}) is the first one with the same assumption that $\alpha$ is helix. The possible solutions to that has already been discussed in Corollary \ref{sonuc1}.
\end{corollary}

\begin{theorem}
Let $\alpha^*$ be the $N-R^*$ associated curve of $\alpha$. The curvature, $\kappa^*$ and the torsion, $\tau^*$ of $\alpha^*$ are given as follows.
\begin{equation}
\begin{aligned}
\kappa^*&=\frac{e^3 \mathbf{L} }{f(e^2+f^2) (\lambda')^3},\\
\tau^*&= \frac{\mathbf{L}^2 (e^2+f^2)}{f^2} \bigg( \mathbf{K} \left( \lambda \kappa^{3}+\lambda \kappa  \tau^{2}-\lambda  \kappa'' - \kappa^{2}-3\, \lambda'' \kappa  -3\, \lambda' \kappa'   \right) \\
& \qquad \qquad \qquad +\mathbf{L } \left( -3\,\lambda  \kappa  \kappa'  -3\,\lambda  \tau  \tau'  -3\, \kappa^{2}\lambda'  -3\, \lambda' \tau^{2}+\lambda'''  +\kappa'   \right) \\
&\qquad \qquad \qquad +\mathbf{M  } \left( -\lambda  \kappa^{2} \tau  -\lambda \tau^{3}+\lambda  \tau''  +\kappa \tau  +3\, \lambda'' \tau  +3\, \lambda' \tau'   \right)  \bigg) \nonumber
\end{aligned}
\end{equation}
\end{theorem}

\begin{prof}
By recalling the relations (\ref{baginti1nr}) and the third derivative (\ref{turev3n}) to substitute into curvatures definitions, (\ref{curvatures}), we complete the proof.
\end{prof}

\section{Binormal Associated Curves}
In this section, we define binormal associated curves such that the binormal vector of a given curve lies on the osculating, normal and rectifying plane of its mate.
\theoremstyle{definition}
\begin{definition}
Let $\alpha$ and $\alpha^*$ be two differentiable curves. If the binormal, B of $\alpha$ is linearly dependent with the vector, $O^*$, then we name the curve $\alpha^*$ as $B-O^*$ associated curve of $\alpha$.
\end{definition}

\begin{theorem}
Let $\alpha^*$ be $B-O^*$ associated curve of $\alpha$ and $\{T^*,N^*,B^*\}$ be the Frenet frame of $\alpha^*$. The relationship of the corresponding Frenet frames are given as follows:
\begin{align}
{T^*}&= \frac{1}{\sqrt { {a}^{2}+{b}^{2}}}\left(\frac {b}{\sqrt{ 1+ \lambda^{2} \tau^{2}} }T-\frac {\lambda  \tau  b}{\sqrt{ 1+ \lambda^{2} \tau^{2}}}N+aB \right)\nonumber\\
{N^*} &=\frac{b}{\sqrt { {a}^{2}+{b}^{2}}}\left(-\frac { \lambda' \mathbf{Y}}{a(\lambda \tau \mathbf{X}+\mathbf{Y})}T+\frac {\lambda  \mathbf{X}}{\lambda \tau \mathbf{X}+\mathbf{Y}}N+B \right) \nonumber\\
{B^*}&= -\frac{b\lambda'}{a (\lambda \tau \mathbf{X}+\mathbf{Y})}\left(\mathbf{X} T+ \mathbf{Y} N\right)\nonumber,
\end{align}
where the coefficients $\mathbf{X} \text{ and } \mathbf{Y}$ are
\begin{align}
\mathbf{X}&=-\lambda  \tau   \left( -\lambda \tau^{2}+\lambda'' \right) - \lambda' \left( -\lambda  \tau'  -2\, \lambda' \tau  +\kappa   \right) \nonumber \\
\mathbf{Y}&= \lambda   \tau^{2}-\lambda''  + \lambda' \lambda  \tau  \kappa \nonumber.
\end{align}
\end{theorem}

\begin{prof}
Since $\alpha$ and $\alpha^*$ are defined as $B-O^*$ associated curves, we write
\begin{equation}\label{bagla3}
\alpha^*=\alpha + \lambda B.
\end{equation}
By differentiating the relation (\ref{bagla3}), using the Frenet formulae given in (\ref{freneForms}) and taking the norm, we have:
\begin{equation}\label{hamtyldz3}
T^*= {\frac {T- \lambda  \tau N + \lambda' B}{\sqrt {1+ \lambda^{2} \tau^{2}+(\lambda')^{2}}}}
\end{equation}
Next taking the second derivative of the equation (\ref{bagla3}) and referring again to (\ref{freneForms}) result the following relation.
\begin{equation}
{\alpha^*}'' =  (\lambda  \tau  \kappa)T +(-\lambda  \tau'  -2\, \lambda' \tau  +\kappa)N +(-\lambda \tau^{2}+\lambda'')B \nonumber.
\end{equation}
The cross production of ${\alpha^*}'$ and ${\alpha^*}''$ leads us the following form,
\begin{equation}
{\alpha^*}' \wedge \, \, {\alpha^*}''= \mathbf{X} T+\mathbf{Y} N+\mathbf{Z}B \nonumber
\end{equation}
where  $\mathbf{X}, \quad \mathbf{Y} \quad \text{and} \quad \mathbf{Z}$ are assigned to be as
\begin{align}\label{kisaltmalar2}
\mathbf{X}&=-\lambda  \tau   \left( -\lambda \tau^{2}+\lambda'' \right) - \lambda' \left( -\lambda  \tau'  -2\, \lambda' \tau  +\kappa   \right) \nonumber \\
\mathbf{Y}&= \lambda   \tau^{2}-\lambda''  + \lambda' \lambda  \tau  \kappa\\
\mathbf{Z}&= -\lambda \tau'  -2\, \lambda' \tau  +\kappa  + \lambda^{2} \tau^{2}\kappa\nonumber
\end{align}
for the sake of simplicity. Note that the norm, $\parallel {\alpha^*}' \wedge \, \, {\alpha^*}''\parallel=\sqrt{\mathbf{X}^2+\mathbf{Y}^2+\mathbf{Z}^2}$.\\ 
By referring again the definitions given by (\ref{defns}), we simply calculate $N^*$, and $B^*$ as
\begin{align}\label{hamnbyldz3}
N^* &= \frac{ (\mathbf{Y} \lambda'+\mathbf{Z} \lambda \tau)T+(-\mathbf{X} \lambda  +\mathbf{Z} )N+(-\mathbf{X}  \lambda  \tau  -\mathbf{Y})B}{\sqrt{1+ \lambda^{2} \tau^{2}+(\lambda')^{2}} \sqrt{\mathbf{X}^2+\mathbf{Y}^2+\mathbf{Z}^2}}\\
B^* &= \frac{\mathbf{X}  T+ \mathbf{Y} N + \mathbf{Z} B}{\sqrt{\mathbf{X}^2+\mathbf{Y}^2+\mathbf{Z}^2}}\nonumber
\end{align}
The intuitive idea is as same as before. Since we defined $\alpha^*$ to be as the $B-O^*$ associated curve of $\alpha $ we can write that
\begin{itemize}
\item $<B,T^*>\,=\,<O^*,T^*>$,
\item $<B,N^*>\,=\,<O^*,N^*> $,
\item $<B,B^*>\,=\,<O^*,B^*>=0$.
\end{itemize}
By using these together with the relations (\ref{oyldz}) and (\ref{hamtyldz3}) we write
\begin{align}\label{baginti1bo}
\bullet \qquad &{\frac {\lambda'}{\sqrt {1+ \lambda^{2} \tau^{2}+(\lambda')^{2}}}}=\frac{a}{\sqrt{a^2+b^2}}, \nonumber \\
\bullet \qquad &\frac{-\mathbf{X}  \lambda  \tau  -\mathbf{Y}}{\sqrt{1+ \lambda^{2} \tau^{2}+(\lambda')^{2}} \sqrt{\mathbf{X}^2+\mathbf{Y}^2+\mathbf{Z}^2}}=\frac{b}{\sqrt{a^2+b^2}}, \\
\bullet \qquad &\frac{\mathbf{Z}}{\sqrt{\mathbf{X}^2+\mathbf{Y}^2+\mathbf{Z}^2}}=0.\nonumber
\end{align}
Substituting these relations into (\ref{hamtyldz3}) and (\ref{hamnbyldz3}), we complete the proof.
\end{prof}
\begin{corollary}\label{sonuc2}
If $\tau$ is taken to be constant then the first relation given in (\ref{baginti1bo}) have the following form
\[\lambda'(\lambda''-\lambda \frac{a^2}{b^2} \tau^2)=0.\]
This relation holds either $\lambda'=0$, correspondingly that $\lambda$ is constant or
\[\lambda=c_1 e^{\frac{a\tau}{b}}+c_2 e^{-\frac{a\tau}{b}}\]
as a result of the solution of second order differential equation, where $c_1$ and $c_2$ are the integration constants. If $\lambda$ is taken to be constant then by the first relation of (\ref{baginti1bo}) $a=0$, resulting that $O^*=N^*$. We remind that this is the definition of Mannheim curves.\\
On the other hand, when considered the third equation in (\ref{baginti1bo}) and recall (\ref{kisaltmalar2}), we have the following
\[\mathbf{Z}=-\lambda \tau'  -2\, \lambda' \tau  +\kappa  + \lambda^{2} \tau^{2}\kappa=0.\]
Rearranging this equation by dividing each term with ($-2\tau$) results
\begin{equation}\label{riccati}
\lambda' +\lambda^2 \left(\frac{-\tau \kappa}{2}\right)+\lambda \left(\frac{\tau'}{2\tau}\right)-\frac{\kappa}{2\tau} =0
\end{equation}
which is clearly a Riccati type of differential equation. If $\lambda=\lambda_1$ is a particular solution for (\ref{riccati}) then we have a general solution by substituting $\lambda=\lambda_1+\frac{1}{\mu}$, that converts the Riccati equation into the fallowing first order linear differential equation:
\begin{equation}\label{convertRiccatitoLinear}
\mu' - \left(2 \lambda_1 \left(\frac{-\tau \kappa}{2}\right)+\left(\frac{\tau'}{2\tau}\right)\right)\mu=\left(\frac{-\tau \kappa}{2}\right)
\end{equation}
where $\mu$ is an arbitrary function of the parameter, $s$. The solution for this (\ref{convertRiccatitoLinear}) can be done by following the steps given in the proof of Theorem (\ref{solODE}).
\end{corollary}

\begin{theorem}
Let $\alpha^*$ be the $B-O^*$ associated curve of $\alpha$. The curvature, $\kappa^*$ and the torsion, $\tau^*$ of $\alpha^*$ are given as follows.
\begin{equation}
\begin{aligned}
\kappa^*&=-\frac{a^4(\mathbf{X}  \lambda  \tau +  \mathbf{Y})}{(\lambda')^2b(a^2+b^2)\sqrt{a^2+b^2}},\\
\tau^*&=\frac{b^2 (\lambda')^2}{a^2 (\mathbf{X}  \lambda  \tau +  \mathbf{Y})^2} \bigg(  \mathbf{X} \left( \lambda  \tau  \kappa'  +3\, \lambda' \tau  \kappa  +2\,\lambda  \tau' \kappa  - \kappa^{2} \right) \\
& \qquad \qquad +\mathbf{Y} \left( \lambda \tau^{3}+\lambda  \tau  \kappa^{2}-\lambda  \tau''  -3\, \lambda' \tau'  -3\, \lambda'' \tau  +\kappa'   \right)\bigg).\nonumber
\end{aligned}
\end{equation}
\end{theorem}
\begin{prof}
By taking the third derivative of (\ref{bagla3}) and using Frenet formulas, we have 
\begin{align}\label{turev3b}
{\alpha^*}'''&=(\lambda \tau  \kappa'  +2\,\lambda \tau' \kappa  +3\, \lambda' \tau \kappa -\kappa^{2})T+(\lambda  \tau \kappa^{2}+\lambda \tau^{3}-\lambda  \tau''  -3\, \lambda'' \tau  -3\, \lambda' \tau'  +\kappa')N \nonumber \\ 
&\qquad +(\lambda'''-3\,\lambda  \tau  \tau'  -3\, \lambda' \tau^{2}+\kappa  \tau)B.
\end{align}
Now, using the relations given in (\ref{baginti1bo}) together with (\ref{turev3b}), to substitute into the definitions (\ref{curvatures}) lets us to complete the proof.
\end{prof}
\theoremstyle{definition}
\begin{definition}
Let $\alpha$ and $\alpha^*$ be two differentiable curves. If the binormal, B of $\alpha$ is linearly dependent with the vector, $P^*$, then we name the curve $\alpha^*$ as $B-P^*$ associated curve of $\alpha$.
\end{definition}

\begin{theorem}
Let $\alpha^*$ be $B-P^*$ associated curve of $\alpha$ and $\{T^*,N^*,B^*\}$ be the Frenet frame of $\alpha^*$. The relationship of the corresponding Frenet Frames are given as follows:
\begin{align}
{T^*}&= \frac {1}{\sqrt {1+{\lambda}^{2} \tau^{2}}}T-\frac {\lambda\,\tau  }{\sqrt {1+{\lambda}^{2} \tau^{2}}}N \nonumber\\
{N^*} &= \frac{1}{\sqrt { {c}^{2}+{d}^{2}}}\left(\frac {d\lambda \tau}{\sqrt {1+{\lambda}^{2} \tau^{2}}}T-\frac {d(\lambda^3 \tau^3-(-\lambda \tau' +\kappa  + \lambda^{2} \tau^{2}))}{(-\lambda \tau' +\kappa  + \lambda^{2} \tau^{2}) \sqrt {1+{\lambda}^{2} \tau^{2}}}N+cB \right) \nonumber\\
{B^*} &= \frac{1}{\sqrt { {c}^{2}+{d}^{2}}}\left(-\frac {c\lambda \tau}{\sqrt {1+{\lambda}^{2} \tau^{2}}}T-\frac{c}{\sqrt {1+{\lambda}^{2} \tau^{2}}}N+dB \right).\nonumber
\end{align}
\end{theorem}

\begin{prof}
Now, since again we defined $\alpha^*$ to be as the $B-P^*$ associated curve of $\alpha $ we can write three of our associative relations as usual which are
\begin{itemize}
\item $ <B,N^*>\,=\,<P^*,N^*>$,
\item $ <B,B^*>\,=\,<P^*,B^*>$,
\item $ <B,T^*>\,=\,<P^*,T^*> =0$.
\end{itemize}
These relations this time result the following three equations
\begin{align}\label{baginti1bp}
\bullet \qquad &\frac{-\mathbf{X}  \lambda  \tau  -\mathbf{Y}}{\sqrt{1+ \lambda^{2} \tau^{2}+(\lambda')^{2}} \sqrt{\mathbf{X}^2+\mathbf{Y}^2+\mathbf{Z}^2}}=\frac{c}{\sqrt{c^2+d^2}}, \nonumber \\
\bullet \qquad &\frac{\mathbf{Z}}{\sqrt{\mathbf{X}^2+\mathbf{Y}^2+\mathbf{Z}^2}}=\frac{d}{\sqrt{c^2+d^2}}, \\
\bullet \qquad &{\frac {\lambda'}{\sqrt {1+ \lambda^{2} \tau^{2}+(\lambda')^{2}}}}=0.\nonumber
\end{align}
When substituted these relations, (\ref{baginti1bp}) in  (\ref{hamtyldz3}) and (\ref{hamnbyldz3}), we complete the proof.
\end{prof}
\begin{corollary}
When taken into account the third relation of (\ref{baginti1bp}) we conclude that if the binormal vector of a given curve is linearly dependent with the unit vector lying on the normal plane of its mate, then the distance between these curves is constant.
\end{corollary}
\begin{theorem}
Let $\alpha^*$ be the $B-P^*$ associated curve of $\alpha$. The curvature, $\kappa^*$ and the torsion, $\tau^*$ of $\alpha^*$ are given as follows.
\begin{equation}
\begin{aligned}
\kappa^*&= -\frac{d^2\sqrt{c^2+d^2} \left(\lambda \tau^2(1+ \lambda^{2} \tau^{2})\right)^3}{c^3 (-\lambda \tau' +\kappa  + \lambda^{2} \tau^{2})^2},\\
\tau*&= \frac{d^2(\lambda^2 \tau^3 \left( \lambda\,\tau  \kappa'  +2\, \tau' \kappa  \lambda- \kappa^{2} \right) + \lambda  \tau^2 \left( \lambda\,\tau   \kappa^{2}+\lambda\, \tau^{3}- \tau'' \lambda+\kappa'   \right) +(-\lambda \tau' +\kappa  + \lambda^{2} \tau^{2})\left( -3\, \tau' \tau \lambda+\kappa  \tau   \right))}{(c^2+d^2)(-\lambda \tau' +\kappa  + \lambda^{2} \tau^{2})^2 } \nonumber
\end{aligned}
\end{equation}
 
\end{theorem}

\begin{prof}
By substituting the relations given in (\ref{baginti1nr}) and the third derivative (\ref{turev3b}) into the definitions given in (\ref{curvatures}), we complete the proof.
\end{prof}

\theoremstyle{definition}
\begin{definition}
Let $\alpha$ and $\alpha^*$ be two differentiable curves. If the binormal, B of $\alpha$ is linearly dependent with the vector, $R^*$, then we name the curve $\alpha^*$ as $B-R^*$ associated curve of $\alpha$.
\end{definition}

\begin{theorem}
Let $\alpha^*$ be $B-R^*$ associated curve of $\alpha$ and $\{T^*,N^*,B^*\}$ be the Frenet frame of $\alpha^*$. The relationship of the corresponding Frenet Frames are given as follows:
\begin{equation}
\begin{aligned}
{T^*}&= \frac{1}{\sqrt { {e}^{2}+{f}^{2}}}\left(\frac {f}{\sqrt{ 1+ \lambda^{2} \tau^{2}} }T-\frac {\lambda  \tau  f}{\sqrt{ 1+ \lambda^{2} \tau^{2}}}N+e B \right), \nonumber \\
{N^*} &=\frac{ef}{\sqrt{e^2+f^2}}\bigg( \frac{\mathbf{Y}\lambda' + \mathbf{Z} \lambda \tau}{\mathbf{Z}\lambda'}T+\frac{-\mathbf{X}\lambda + \mathbf{Z}}{\mathbf{Z} \lambda'} \bigg), \nonumber \\
{B^*} &= \frac{f}{\sqrt{e^2+f^2}}\bigg( \frac{\mathbf{X}}{\mathbf{Z}}T +\frac{\mathbf{Y}}{\mathbf{Z}}N + B \bigg).
\end{aligned}
\end{equation}

\end{theorem}
\begin{prof}
Now, since again we defined $\alpha^*$ to be as the $B-R^*$ associated curve of $\alpha $ we can write three of our associative relations as usual which are
\begin{itemize}
\item $ <B,T^*>\,=\,<R^*,T^*>$,
\item $ <B,B^*>\,=\,<R^*,B^*>$,
\item $ <B,N^*>\,=\,<R^*,N^*> =0$.
\end{itemize}
These relations this time result the following three equations
\begin{align}\label{baginti1br}
\bullet \qquad &{\frac {\lambda'}{\sqrt {1+ \lambda^{2} \tau^{2}+(\lambda')^{2}}}}=\frac{e}{\sqrt{e^2+f^2}}, \nonumber \\
\bullet \qquad &\frac{\mathbf{Z}}{\sqrt{\mathbf{X}^2+\mathbf{Y}^2+\mathbf{Z}^2}}=\frac{f}{\sqrt{e^2+f^2}}, \\
\bullet \qquad &\frac{-\mathbf{X}  \lambda  \tau  -\mathbf{Y}}{\sqrt{1+ \lambda^{2} \tau^{2}+(\lambda')^{2}} \sqrt{\mathbf{X}^2+\mathbf{Y}^2+\mathbf{Z}^2}}=0.\nonumber
\end{align}
For the last time when substituted (\ref{baginti1br}) into (\ref{hamtyldz3}) and (\ref{hamnbyldz3}), the proof is complete.
\end{prof}
\begin{corollary}
The only analytically solvable equation is the first one of (\ref{baginti1br}) with the same assumption given in Corollary (\ref{sonuc2}). The possible solutions can be get by fallowing the same steps as well.
\end{corollary}
\begin{theorem}
Let $\alpha^*$ be the $B-R^*$ associated curve of $\alpha$. The curvature, $\kappa^*$ and the torsion, $\tau^*$ of $\alpha^*$ are given as follows.
\begin{equation}
\begin{aligned}
\kappa^*&=\frac{\mathbf{Z}e^3}{f(e^2+f^2)(\lambda')^3},\\
\tau^*&=\frac{f^2}{\mathbf{Z}^2 (e^2+f^2)}\bigg( \mathbf{X} \left( \lambda  \tau  \kappa'  +2\,\lambda   \tau' \kappa  +3\,\lambda' \tau  \kappa  - \kappa^{2} \right) +\mathbf{Y} \left( \lambda  \tau \kappa^{2}+\lambda \tau^{3}-\lambda \tau''  -3\,\lambda'' \tau  -3\, \lambda' \tau'  +\kappa'   \right) \\
&\qquad \qquad +\mathbf{Z} \left( -3\,\lambda  \tau  \tau'  -3\, \left( \lambda'   \right) \tau^{2}+\kappa  \tau  +\lambda'''   \right) \bigg). \nonumber
\end{aligned}
\end{equation}
\end{theorem}
\begin{prof}
Recall the relations (\ref{baginti1br}) and (\ref{turev3b}), and substitute these in (\ref{curvatures}), the proof is complete.
\end{prof}


\end{document}